\documentclass[twoside,10pt]{article}

\usepackage{amsmath}
\usepackage{amsfonts} 
\usepackage{cite}

\begin{document}

\title{Nonlinear Sequence Transformations: \\ 
Computational Tools for the  Acceleration of Convergence \\ 
and the Summation of Divergent Series}

\author{Ernst Joachim Weniger \\    
Institut f\"ur Physikalische und Theoretische Chemie \\
Universit\"at Regensburg, D-93040 Regensburg, Germany \\
joachim.weniger@chemie.uni-regensburg.de}

\date{Sububmitted to the topical issue \\ 
\textit{Time Series and Time Evolution of Generic Systems. \\
Spectral Analysis beyond Fourier Transform, \\
Auto- and Cross-Correlation Functions, Green's Resolvent Spectra. \\ 
Nonlinear Convergence Accelerators of Series and Sequences} \\
of \\ 
Journal of Computational Methods in Sciences and Engineering (JCMSE) 
\\ Corrected Version 11 June 2001}

\maketitle

\typeout{==> Abstract}
\begin{abstract}
\noindent
Convergence problems occur abundantly in all branches of mathematics or
in the mathematical treatment of the sciences. Sequence transformations
are principal tools to overcome convergence problems of the kind. They
accomplish this by converting a slowly converging or diverging input
sequence $\{ s_n \}_{n=0}^{\infty}$ into another sequence $\{
s^{\prime}_n \}_{n=0}^{\infty}$ with hopefully better numerical
properties. Pad\'{e} approximants, which convert the partial sums of a
power series to a doubly indexed sequence of rational functions, are the 
best known sequence transformations, but the emphasis of the review will 
be on alternative sequence transformations which for some problems
provide better results than Pad\'{e} approximants.
\end{abstract} 

\typeout{==> Section 1}
\section{Introduction}
\label{Sec:Intro}

Many numerical techniques as for example iterative schemes,
discretization methods, perturbation techniques, or series expansions
produce results which are actually \emph{sequences}. Obviously, a
numerical technique of that kind is practically useful only if the
resulting sequence converges sufficiently fast. Unfortunately, it
frequently happens that the resulting sequence either converges
\emph{too slowly} to be practically useful, or it may even
\emph{diverge}.

Problems with slow convergence or divergence were of course already
encountered in the early days of calculus. Accordingly, numerical
techniques for the acceleration of convergence or the summation of
divergent series are almost as old as calculus itself. According to
Knopp \cite[p. 249]{Kno1964}, the first series transformation was
published by Stirling \cite{Sti1730} already in 1730, and in 1755 Euler
\cite{Eul1755} published the series transformation which now bears his
name. In rudimentary form, convergence acceleration methods are even
older. In a book by Brezinski \cite[pp.\ 90 - 91]{Bre1991a} it is
mentioned that convergence acceleration methods were already used in
1654 by Huygens and in 1674 by Seki Kowa, the probably most famous
Japanese mathematician of his time. Both Huygens and Seki Kowa tried to
obtain better approximations to $\pi$. Huygens used a linear
extrapolation scheme which is a special case of what we now call
Richardson extrapolation \cite{Ri1927}, and Seki Kowa used the so-called
$\Delta^2$ process, which is usually attributed to Aitken
\cite{Ai1926}. Then, in a book by Liem, L\"{u}, and Shih \cite[p.\
$ix$]{LieLueSh1995} it is mentioned that extrapolation methods were
already used by the Chinese mathematicians Liu Hui (A.D. 263) and Zhu
Chongzhi (429 - 500) for obtaining better approximations to $\pi$, but
no further details are given.

\emph{Sequence transformations} are principal tools to overcome
convergence problems of the kind mentioned above. In this approach, a
slowly convergent or divergent sequence $\{ s_n \}_{n=0}^{\infty}$,
whose elements may for instance be the partial sums
\begin{equation}
s_n \; = \; \sum_{k=0}^{n} \, a_k
\label{Par_Sum}
\end{equation}
of an infinite series, is converted into a new sequence $\{ s^{\prime}_n
\}_{n=0}^{\infty}$ with hopefully better numerical properties. The
history of sequence transformations and related topics starting from the
17th century until today is discussed by Brezinski in a monograph
\cite{Bre1991a} or in two articles \cite{Bre1996,Bre2000}.

Before the invention of electronic computers, mainly \emph{linear}
sequence transformations were used, which compute the elements of the
transformed sequence $\{ s_n^{\prime} \}_{n=0}^{\infty}$ as weighted
averages of the elements of the input sequence $\{ s_n
\}_{n=0}^{\infty}$ according to
\begin{equation}
s_n^{\prime} \; = \; \sum_{k=0}^{n} \, \mu_{n k} \, s_k \, .
\label{Mat_tr}
\end{equation}
The theoretical properties of these matrix transformations are now very
well understood
\cite{Kno1964,Ha1949,Pet1966,Pey1969,ZelBee1970,PoSh1988,Wi1981}. Their 
main appeal lies in the fact that for the weights ${\mu}_{n k}$ in
(\ref{Mat_tr}) some necessary and sufficient conditions could be
formulated which guarantee that the application of such a matrix
transformation to a \emph{convergent} sequence $\{ s_n
\}_{n=0}^{\infty}$ yields a transformed sequence $\{ s_n^{\prime}
\}_{n=0}^{\infty}$ converging to the \emph{same} limit $s =
s_{\infty}$. Theoretically, this \emph{regularity} is extremely
desirable, but from a practical point of view, it is a
disadvantage. This sounds paradoxical. However, Wimp remarks in the
preface of his book \cite[p.\ X]{Wi1981} that the size of the domain of
regularity of a transformation and its efficiency seem to be inversely
related. Accordingly, regular linear transformations are in general at
most moderately powerful, and the popularity of most linear
transformations has declined considerably in recent years.

Modern \emph{nonlinear} sequence transformations as for instance Wynn's
epsilon \cite{Wy1956a} and rho \cite{Wy1956b} algorithm or Brezinski's
theta algorithm \cite{Bre1971} have largely complementary properties:
They are nonregular, which means that that the convergence of the
transformed sequence to the correct limit is not guaranteed. In
addition, their theoretical properties are far from being completely
understood. Nevertheless, they often accomplish spectacular
results. Consequently, nonlinear transformations now dominate both
mathematical research as well as practical applications, as documented
by the large number of recent books
\cite{LieLueSh1995,Wi1981,Bre1977,Bre1978,Bre1980,MarSha1983,%
CuWu1987,De1988,BreRZa1991,Wal1996} and review articles
\cite{Gu1989,We1989,Hom2000} on this topic.

There is considerable evidence that the culprit for the frequently
unsatisfactory performance of regular matrix transformation is not their
linearity, but their regularity. In \cite{We1994a} it was shown that
suitably chosen linear but nonregular transformations can at least for
special problems be as efficient as the most powerful nonlinear
transformations.

The best known class of sequence transformations are Pad\'{e}
approximants which convert the partial sums
\begin{equation}
f_n (z) \; = \; \sum_{k=0}^{n}\, \gamma_k \, z^k
\label{Par_Sum_PS}
\end{equation}
of a (formal) power series for some function $f$ into a doubly indexed
sequence of rational functions
\begin{equation}
[ l / m ]_f (z) \; = \; P_{l} (z) / Q_m (z) \, , 
\qquad l, m \in \mathbb{N}_0 \, .
\end{equation}
Here, $P_{l} (z) = p_0 + p_1 z + \ldots + p_{l} z^{l}$ and $Q_m (z)= 1 +
q_1 z + \ldots + q_m z^m$ are polynomials in $z$ of degrees $l$ and $m$,
respectively. The $l + m + 1$ polynomial coefficients $p_0$, $p_1$,
$\ldots$, $p_{l}$ and $q_1$, $q_2$, $\ldots$, $q_m$ are chosen in such a
way that the Taylor expansion of the ratio $P_{l} (z) / Q_m (z)$ at $z =
0$ agrees with the power series for $f$ as far as possible:
\begin{equation}
f (z) \, - \, P_{l} (z) / Q_m (z) \; = \;
O (z^{l + m +1}) \, , \quad z \to 0 \, .
\label{O_est_Pade}
\end{equation}
This asymptotic error estimate leads to a system of $l + m + 1$ linear
equations for the coefficients of the polynomials $P_{l} (z)$ and $Q_m
(z)$ \cite{Ba1975,BaGM1996}. Moreover, several recursive algorithms are
known. The merits and weaknesses of the various computational schemes
for Pad\'{e} approximants are discussed in \cite[Section
II.3]{CuWu1987}.

In applied mathematics and in theoretical physics, Pad\'{e} approximants
have become the standard tool to overcome convergence problems with
power series. Accordingly, there is an extensive literature, and any
attempt to provide a complete bibliography would be beyond the scope of
this article (see for example the extensive bibliography compiled by
Brezinski \cite{Bre1991b}). The popularity of Pad\'{e} approximants in
theoretical physics can be traced back to an article by Baker
\cite{Ba1965}, who also wrote the first modern monograph \cite{Ba1975}. 
The currently most complete treatment of the theory of Pad\'{e}
approximants is the 2nd edition of the monograph by Baker and
Graves-Morris \cite{BaGM1996}. More condensed treatments can be found in
books on continued fractions
\cite{Wal1973,JonThr1980,BowShe1989,LorWaa1992}, in a book by Bender and
Orszag on mathematical physics \cite[Section 8]{BenOrs1978}, or in a
book by Baker on critical phenomena \cite[Part III]{Ba1990}. Then, there
is a book by Pozzi on the use of Pad\'{e} approximants in fluid dynamics
\cite{Po1994} as well as several review articles
\cite{ZiJ1971,Bas1972,BrIs1994,BrIs1995}. The generalization of Pad\'{e} 
approximants to operators and multivariate power series is discussed in
a book by Cuyt \cite{Cuy1984} and in articles by Cuyt \cite{Cuy1999} and
by Guillaume and A.\ Huard \cite{GuiHua2000}. Another generalization of
Pad\'{e} approximants -- Pad\'{e}-type approximants -- are described in
a book by Brezinski \cite{Bre1980}.

The emphasis of this article is not on Pad\'{e} approximants, which are
well known as well as extensively documented in the literature, but on
alternative sequence transformations, which are not so well known
yet. This is quite undeserved. In some cases, sequence transformations
outperform Pad\'{e} approximants. For example, the present author has
applied sequence transformations successfully in such diverse fields as
the evaluation of special functions
\cite{We1989,We1994a,We1990,WeCi1990,We1996c,JeMoSoWe1999,We2001}, the
evaluation of molecular multicenter integrals of exponentially decaying
functions \cite{WeGroSte1986,GroWeSte1986,WeSte1989,SteWe1990,HoWe1995},
the summation of strongly divergent quantum mechanical perturbation
expansions \cite{We1990,WeCiVi1991,CiViWe1991,We1992,WeCiVi1993,%
CiViWe1993,We1996a,We1996b,CiWeBraSpi1996,We1996d,We1997,JeBeWeSo2000,%
JeWeS02000}, and the extrapolation of crystal orbital and cluster
calculations for oligomers to their infinite chain limits of
stereoregular quasi-onedimensional organic polymers
\cite{WeLie1990,CioWe1993,WeKi2000}. In the majority of these
applications, it was either not possible to use Pad\'{e} approximants,
or sequence transformations did a better job.

It is of course impossible to present something as complex and diverse
as sequence transformations in a single and comparatively short
article. Because of obvious space limitations, this article can at best
provide some basic facts about some of the most useful sequence
transformations and review the recent literature. However, it cannot be
a substitute for more detailed treatments as for example an older review
of the author \cite{We1989} or the book by Brezinski and Redivo Zaglia
\cite{BreRZa1991}.

\typeout{==> Section 2}
\setcounter{equation}{0}
\section{On the Construction of Sequence Transformations}
\label{Sec:ConstrSeqTran}

The basic step for the construction of a sequence transformation is the
assumption that the elements of a convergent or divergent sequence $\{
s_n \}_{n=0}^{\infty}$ can be partitioned into a \emph{limit} or
\emph{antilimit} $s$ and a \emph{remainder} $r_n$ according to
\begin{equation}
s_n \; = \; s + r_n \, , \qquad n \in \mathbb{N}_0 \, .
\label{s_n_r_n} 
\end{equation}
If a sequence $\{ s_n \}_{n=0}^{\infty}$ converges, the remainders $r_n$
in (\ref{s_n_r_n}) can be made negligible by increasing $n$ as much as
necessary. However, many sequences converge so slowly that this is not
feasible. Moreover, increasing $n$ does not help in the case of
divergence.

Alternatively, one can try to improve convergence by computing
approximations to the remainders which are then eliminated from the
sequence elements. Of course, this is more easily said than done. The
remainders $\{ r_n \}_{n=0}^{\infty}$ of a sequence $\{ s_n
\}_{n=0}^{\infty}$ are in general unknown, and their determination is
normally not easier than the determination of the (generalized) limit.
For example, if the input data $\{ s_n \}_{n=0}^{\infty}$ are the
partial sums (\ref{Par_Sum}) of an infinite series, the remainders
satisfy 
\begin{equation}
r_n \; = \; - \sum_{k=n+1}^{\infty} \, a_k \, .
\end{equation}
Thus, the straightforward elimination of exact remainders is normally
not possible. However, the remainders of some infinite series can be
approximated with the help of the Euler-Maclaurin formula (see for
example \cite[pp.\ 279 - 295]{Olv1974}). Let us consider a convergent
infinite series $\sum_{\nu=0}^{\infty} g (\nu)$, and let us assume that
its terms $g (\nu)$ are \emph{smooth} and \emph{slowly varying}
functions of the index $\nu$. Then, the integral
\begin{equation}
\int_{M}^{N} \, g (x) \, \mathrm{d} x
\label{EuMacInt}
\end{equation}
with $M, N \in \mathbb{Z}$ provides a good approximation to the sum
\begin{equation}
\frac{1}{2} g (M) \, + \, g (M+1) \, + \, \ldots \, + \, g (N-1) \, + \,
\frac{1}{2} g (N)
\end{equation}
and vice versa. In the years between 1730 and 1740, Euler and Maclaurin
derived independently correction terms, which ultimately yielded what we
now call the Euler-Maclaurin formula:
\begin{subequations}
\label{EuMaclau}
\begin{align}
\sum_{\nu=M}^{N} \, g (\nu) & \; = \;
\int_{M}^{N} \, g (x) \, \mathrm{d} x
\, + \, \frac{1}{2} \left[ g (M) + g (N) \right] \notag \\
& \qquad \, + \, \sum_{j=1}^{k} \, \frac {B_{2j}} {(2j)!}
\left[ g^{(2j-1)} (N) - g^{(2j-1)} (M) \right] \, + \, R_k (g) \, , \\
R_k (g) & \; = \; - \, \frac{1}{(2k)!} \,
\int_{M}^{N} \, B_{2k}
\bigl( x - [\mkern - 2.5 mu  [x] \mkern - 2.5 mu ] \bigr) \,
g^{(2k)} (x) \, \mathrm{d} x \, .
\end{align}
\end{subequations}
Here, $[\mkern - 2.5 mu [x] \mkern - 2.5 mu ]$ is the integral part of
$x$, $B_m (x)$ is a Bernoulli polynomial , and $B_m = B_m (0)$ is a
Bernoulli number.                          

If we set $M = n+1$ and $N = \infty$, the leading terms of the
Euler-Maclaurin formula, which is actually an asymptotic series, yield
for sufficiently large $n$ rapidly convergent approximations to the
truncation error $\sum_{\nu=n+1}^{\infty} g (\nu)$. 

An example, which demonstrates the usefulness of the Euler-Maclaurin
formula, is the Dirichlet series for the Riemann zeta function:
\begin{equation}
\zeta (z) \; = \; \sum_{\nu = 0}^{\infty} \, (\nu + 1)^{-z} \, .
\label{ZetaSer}
\end{equation}
This series converges for $\mathrm{Re} (z) > 1$. However, it is
notorious for extremely slow convergence if $\mathrm{Re} (z)$ is only
slightly larger than one.

The terms $(\nu+1)^{-z}$ of the Dirichlet series (\ref{ZetaSer}) are
obviously smooth and slowly varying functions of the index $\nu$ and
they can be differentiated and integrated easily. Thus, we can apply the
Euler-Maclaurin formula (\ref{EuMaclau}) with $M = n+1$ and $N = \infty$
to the truncation error of the Dirichlet series:
\begin{subequations}
\label{EuMacZeta}
\begin{align}
\sum_{\nu=n+1}^{\infty} \, (\nu+1)^{-z} & \; = \;
\frac{(n+2)^{1-z}}{z-1}
\, + \, \frac{1}{2} \, (n+2)^{-z} \notag \\
& \qquad \, + \,
\sum_{j=1}^{k} \, \frac {(z)_{2j-1} \, B_{2j}} {(2j)!} \,
(n+2)^{-z-2j+1} \, + \, R_k (n,z) \, , \\
R_k (n,z) & \; = \; - \, \frac{(z)_{2k}}{(2k)!} \,
\int_{n+1}^{\infty} \, \frac
{B_{2k} \bigl( x - [\mkern - 2.5 mu  [x] \mkern - 2.5 mu ] \bigr)}
{(1+x)^{z+2k}} \, \mathrm{d} x \, .
\end{align}
\end{subequations}
Here, $(z)_m = z (z+1) \cdots (z+m-1) = \Gamma(z+m)/\Gamma(z)$ is a
Pochhammer symbol.  

In \cite[Tables 8.7 and 8.8, p.\ 380]{BenOrs1978} or in \cite[Section
2]{WeKi2000} it was shown that a few terms of the Euler-Maclaurin
approximation (\ref{EuMacZeta}) suffice for a convenient and reliable
computation of $\zeta (z)$ even if $z$ is so close to one. The arguments
$z = 1.1$ and $z = 1.01$ considered in \cite{BenOrs1978,WeKi2000} lead
to such a slow convergence of the Dirichlet series (\ref{ZetaSer}) that
its evaluation via a straightforward addition of its terms is
practically impossible.

Thus, it looks like an obvious idea to use the Euler-Maclaurin formula
routinely in the case of slowly convergent series. Unfortunately, this
is not possible. The Euler-Maclaurin formula requires that the terms of
the series can be differentiated and integrated with respect to the
index. This excludes many series of interest. Moreover, the
Euler-Maclaurin formula cannot be applied in the case of alternating or
divergent series since their terms are neither smooth nor slowly
varying. However, the probably worst drawback of the Euler-Maclaurin
formula is that it is an \emph{analytic} convergence acceleration
method. This means that it cannot be applied if only the numerical
values of the terms of a series are known.

Sequence transformations also try to compute approximations to the
remainders and to eliminate them from the sequence elements. However,
they require only relatively little knowledge about the $n$-dependence
of the remainders of the sequence to be transformed. Consequently, they
can be applied in situations in which apart from the numerical values of
a finite string of sequence elements virtually nothing else is known.

Since it would be futile to try to eliminate a remainder $r_n$ with a
completely unknown and arbitrary $n$-dependence, a sequence
transformation has to make some assumptions, either implicitly or
explicitly. Therefore, a sequence transformation will only work well if
the actual behavior of the remainders is in sufficient agreement with
the assumptions made. Of course, this also implies that a sequence
transformation, which makes certain assumptions, may fail to accomplish
something if it is applied to a sequence with remainders of a
sufficiently different behavior. Thus, the efficiency of a sequence
transformation for certain sequences and its inefficiency or even
nonregularity for other sequences are intimately related.

Sequence transformations normally eliminate only approximations to the
remainders. In such a case, the elements of the transformed sequence $\{
s_n^{\prime} \}_{n=0}^{\infty}$ will also be of the type of
(\ref{s_n_r_n}), which means that $s^{\prime}_n$ can also be partitioned
into the (generalized) limit $s$ and a transformed remainder
$r_n^{\prime}$ according to
\begin{equation}
s_n^{\prime} \; = \; s + r_n^{\prime} \, , 
\qquad n \in \mathbb{N}_0 \, .
\end{equation}
The transformed remainders $\{ r_n^{\prime} \}_{n=0}^{\infty}$ are
normally different from zero for all finite values of $n$. However,
convergence is \emph{accelerated} if the transformed remainders $\{
r_n^{\prime} \}_{n=0}^{\infty}$ vanish more rapidly than the original
remainders $\{ r_n \}_{n=0}^{\infty}$,
\begin{equation}
\lim_{n \to \infty} \, \frac {s^{\prime}_n - s} {s_n - s} \; = \;
\lim_{n \to \infty} \, \frac {r^{\prime}_n} {r_n} \; = \; 0 \, ,
\end{equation}
and a divergent sequence is summed if the transformed remainders $r'_n$
vanish as $n \to \infty$.

Assumptions about the $n$-dependence of the truncation errors can be
incorporated into the transformation process via \emph{model
sequences}. In this approach, a sequence transformation $E_{k}^{(n)}$ is
constructed in such a way that it produces the (generalized) limit
$\tilde{s}$ of a model sequence
\begin{equation}
\tilde{s}_n \; = \; \tilde{s} \, + \, \tilde{r}_n \; = \; \tilde{s}
\, + \, \sum_{j=0}^{k-1} \, {\tilde{c}}_j \, \phi_j (n)  \, ,
\label{Mod_Seq}
\end{equation}
if it is applied to a set of $k+1$ consecutive elements of this model
sequence:
\begin{equation}
E_{k}^{(n)} \; = \; 
E_{k}^{(n)} (\tilde{s}_{n}, \tilde{s}_{n+1}, \ldots , \tilde{s}_{n+k})
\; = \; \tilde{s} \, .
\end{equation}
The $\phi_j (n)$ are assumed to be \emph{known} functions of $n$,
and the ${\tilde{c}}_j$ are unspecified coefficients.   

The elements of this model sequence contain $k+1$ unknowns, the limit
or antilimit $\tilde{s}$ and the $k$ coefficients $\tilde{c}_0$,
$\tilde{c}_1$, $\ldots$, $\tilde{c}_{k-1}$. Since all unknowns occur
\emph{linearly}, Cramer's rule implies that a sequence
transformation, which is able to determine the limit or antilimit
$\tilde{s}$ in (\ref{Mod_Seq}) from the numerical values of $k+1$
sequence elements $\tilde{s}_n$, $\tilde{s}_{n+1}$, $\ldots$,
$\tilde{s}_{n+k}$, is given as the ratio of two determinants \cite[p.\
56]{BreRZa1991}. Since determinantal representations are computationally
unattractive, alternative computational schemes for $E_{k}^{(n)}$ are
highly desirable. A recursive scheme for $E_{k}^{(n)}$ was derived
independently by Schneider \cite{Schn1975}, Brezinski \cite{Bre1980b},
and H{\aa}vie \cite{Haa1979}. A more economical implementation was later
obtained by Ford and Sidi \cite{ForSid1987}.

The sequence transformation $E_{k}^{(n)}$ contains the unspecified
quantities $\phi_j (n)$. The majority of all currently known sequence
transformations can be obtained by specializing the $\phi_j (n)$ (see
for example \cite[pp.\ 57 - 58]{BreRZa1991}). So, from a purely formal
point of view either the determinantal representation for $E_{k}^{(n)}$
or the recursive schemes provide a complete solution to the majority of
all convergence acceleration and summation problems. However, even the
recursive scheme of Ford and Sidi \cite{ForSid1987} is considerably more
complicated than the recursive schemes for other transformations that
can be obtained by specializing the $\phi_j (n)$ in (\ref{Mod_Seq}). So,
is it usually simpler to use instead of $E_{k}^{(n)}$ one of its special
cases. Important is also the following aspect: It is certainly helpful
to know that for arbitrary functions $\phi_j (n)$ the sequence
transformation $E_k^{(n)}$ can be computed, but it is more important to
find out which set $\{ \phi_{j} (n) \}_{j=0}^{\infty}$ produces the best
results for a given sequence $\{ s_n \}_{n=0}^{\infty}$.

If the remainders of the model sequence (\ref{Mod_Seq}) are capable of
producing sufficiently accurate approximations to the remainders of a
sequence $\{ s_n \}_{n=0}^{\infty}$, then the application of the
sequence transformation $E_{k}^{(n)}$ to $k+1$ sequence elements $s_n$,
$s_{n+1}$, $\ldots$ , $s_{n+k}$ should produce a sufficiently accurate
approximation to the (generalized) limit $s$ of the input sequence. This
is usually the case if the truncation errors can be expressed as
\emph{infinite} series in terms of the $\phi_j (n)$. Further details as
well as many examples can be found in \cite{BreRZa1991,We1989}.

Simple asymptotic conditions are used in the literature to classify the
type of convergence of a sequence. For example, many practically
relevant sequences $\{ s_n \}_{n=0}^{\infty}$ converging to some limit
$s$ can be characterized by the asymptotic condition
\begin{equation}
\lim_{n \to \infty} \, \frac{s_{n+1} - s}{s_n - s} \; = \; \rho \, ,
\label{LinLogConv}
\end{equation}
which closely resembles the ratio test in the theory of infinite series.
If $\vert \rho \vert < 1$, the sequence is called \emph{linearly}
convergent, and if $\rho = 1$, it is called \emph{logarithmically}
convergent.  

Typical examples of linearly convergent sequences are the partial sums
of a power series with a nonzero, but finite radius of convergence. In
contrast, the partial sums $\sum_{\nu = 0}^{n} (\nu + 1)^{-z}$ of the
Dirichlet series (\ref{ZetaSer}) for $\zeta (z)$ converge
logarithmically.

\typeout{==> Section 3}
\setcounter{equation}{0}
\section{The Aitken Formula, Wynn's Epsilon Algorithm, and Related
Transformations} 
\label{Sec:AitEps}

It will now be shown how Aitken's $\Delta^2$ formula can be constructed
by assuming that the truncation error consist of a single exponential
term according to
\begin{equation}
s_n \; = \; s \, + \, c \, \lambda^n \, , \qquad c \ne 0,
\quad \vert \lambda \vert \ne 1 \, , \quad n \in \mathbb{N}_0 \, .
\label{AitModSeq}
\end{equation}
For $0 < \vert \lambda \vert < 1$, this sequence converges to its limit
$s$, and for $\vert \lambda \vert > 1$, it diverges away from its
generalized limit or antilimit $s$. Thus, if this sequence converges, it
converges linearly according to (\ref{LinLogConv}).

By considering $s$, $c$, and $\lambda$ in (\ref{AitModSeq}) as unknowns
of the linear system $s_{n+j} = s + c \lambda^{n+j}$ with $j = 0, 1, 2$,
a sequence transformation can be constructed which is able to determine
the (generalized) limit $s$ of the model sequence (\ref{AitModSeq}) from
the numerical values of three consecutive sequence elements $s_n$,
$s_{n+1}$ and $s_{n+2}$. A short calculation shows that
\begin{equation}
\mathcal{A}_{1}^{(n)} \; = \; s_n \, - \,
\frac{[\Delta s_n]^2}{\Delta^2 s_n} \, ,  \qquad n \in \mathbb{N}_0 \, ,
\label{AitFor_1}
\end{equation}
produces the (generalized) limit $s$ of the model sequence
(\ref{AitModSeq}) according to $\mathcal{A}_{1}^{(n)} = s$. Alternative
expressions for $\mathcal{A}_{1}^{(n)}$ are discussed in \cite[Section
5.1]{We1989}. The forward difference operator $\Delta$ in
(\ref{AitFor_1}) is defined according to $\Delta G (n) = G (n+1 - G
(n)$.

The power and practical usefulness of Aitken's $\Delta^2$ formula is
limited since it is designed to eliminate only a single exponential
term. However, the output data ${\cal A}_1^{(n)}$ can be used as input
data in the $\Delta^2$ formula (\ref{AitFor_1}). Hence, the $\Delta^2$
process can be iterated, which leads to the following nonlinear
recursive scheme \cite[Eq.\ (5.1-15)]{We1989}:
\begin{subequations}
\begin{eqnarray}
\lefteqn{} \nonumber \\
{\cal A}_0^{(n)} & = & s_n \, , \qquad n \in \mathbb{N}_0 \, , \\
{\cal A}_{k+1}^{(n)} & = & {\cal A}_{k}^{(n)} -
\frac
{\bigl[\Delta {\cal A}_{k}^{(n)}\bigr]^2}
{\Delta^2 {\cal A}_{k}^{(n)}} \, , \qquad k, n \in \mathbb{N}_0 \, .
\end{eqnarray}
\label{It_Aitken}
\end{subequations}%
In the case of doubly indexed quantities like $\mathcal{A}_{k}^{(n)}$,
it will always be assumed that $\Delta$ only acts on the superscript $n$
but not on the subscript $k$ according to $\Delta \mathcal{A}_{k}^{(n)}
= \mathcal{A}_{k}^{(n+1)} - \mathcal{A}_{k}^{(n)}$.

There is an extensive literature on Aitken's $\Delta^2$ process and its
iteration. Its numerical performance was studied in
\cite{We1989,We2001,SmiFor1982}. It was also discussed in articles by Lubkin
\cite{Lub1952}, Shanks \cite{Sha1955}, Tucker \cite{Tuc1967,Tuc1969},
Clark, Gray, and Adams \cite{ClGrAd1969}, Cordellier \cite{Cor1979},
Hillion \cite{Hil1975}, Jurkat \cite{Jur1983}, Bell and Phillips 
\cite{BelPhi1984}, and Weniger \cite{We1989,We2001,We2000}, or in books
by Baker and Graves-Morris \cite{BaGM1996}, Brezinski
\cite{Bre1977,Bre1978}, Brezinski and Redivo Zaglia \cite{BreRZa1991},
Delahaye \cite{De1988}, Walz \cite{Wal1996}, and Wimp \cite{Wi1981}. A 
multidimensional generalization of Aitken's transformation to vector
sequences was discussed by MacLeod \cite{MacL1986}. Modifications and
generalizations of Aitken's $\Delta^2$ process were proposed by Drummond
\cite{Dru1976}, Jamieson and O'Beirne \cite{JamOBe1978}, Bj{\o}rstad, 
Dahlquist, and Grosse \cite{BjoDahGro1981}, and Sablonniere
\cite{Sab1992}. The iteration of other sequence transformations is 
discussed in \cite{We1991}.

An obvious generalization of the model sequence (\ref{AitModSeq}) would
be the following model sequence which contains $k$ exponential terms:
\begin{equation}
s_n \; = \; s \, + \, \sum_{j=0}^{k-1} \, c_j \, \lambda_{j}^{n} \, ,
\qquad \vert \lambda_{0} \vert > \vert \lambda_{1} \vert
> \ldots > \vert \lambda_{k-1} \vert \, .
\label{EpsModSeq1}
\end{equation}
Although the $\Delta^2$ process (\ref{AitFor_1}) is by construction
exact for the model sequence (\ref{AitModSeq}), its iteration
(\ref{It_Aitken}) is not exact for the model sequence
(\ref{EpsModSeq1}). Instead, this is -- as first shown by Wynn
\cite{Wy1966} and later extended by Sidi \cite{Sid1996} -- true for 
Wynn's epsilon algorithm \cite{Wy1956a}:
\begin{subequations}
\label{eps_al}
\begin{eqnarray}
\epsilon_{-1}^{(n)} & \; = \; & 0 \, ,
\qquad \epsilon_0^{(n)} \, = \, s_n \, ,
\qquad  n \in \mathbb{N}_0 \, , \\
\epsilon_{k+1}^{(n)} & \; = \; & \epsilon_{k-1}^{(n+1)} \, + \,
1 / [\epsilon_{k}^{(n+1)} - \epsilon_{k}^{(n)} ] \, ,
\qquad k, n \in \mathbb{N}_0 \, .
\end{eqnarray}
\end{subequations}
Only the elements $\epsilon_{2k}^{(n)}$ with \emph{even} subscripts
provide approximations to the limit $s$ of the sequence $\{ s_n
\}_{n=0}^{\infty}$ to be transformed. Efficient recursive schemes based 
on the \emph{moving lozenge technique} introduced by Wynn \cite{Wy1965}
are discussed in \cite[Section 4.3]{We1989}.

A short calculation shows $\mathcal{A}_{1}^{(n)} =
\epsilon_{2}^{(n)}$. For $k > 1$, $\mathcal{A}_{k}^{(n)}$ and
$\epsilon_{2k}^{(n)}$ are in general different, but they have similar
properties in convergence acceleration and summation processes.

If the input data $s_n$ for Wynn's epsilon algorithm are the partial
sums (\ref{Par_Sum_PS}) of the (formal) power series for some function
$f (z)$ according to $s_n = f_n (z)$, then it produces Pad\'{e}
approximants according to \cite{Wy1956a}
\begin{equation}
\epsilon_{2 k}^{(n)} \; = \; [ n + k / k ]_f (z) \, , 
\qquad k, n \in \mathbb{N}_0 \, .
\label{Eps_Pade}
\end{equation}

Since the epsilon algorithm can be used for the computation of Pad\'e
approximants, it is discussed in books on Pad\'e approximants such as
\cite{BaGM1996}, but there is also an extensive literature dealing
directly with it. In Wimp's book \cite[p.\ 120]{Wi1981}, it is mentioned
that over 50 articles on the epsilon algorithm were published by Wynn
alone, and at least 30 articles by Brezinski. As a fairly complete
source of references Wimp recommends Brezinski's first book
\cite{Bre1977}. However, this book was published in 1977, and since then
many more articles dealing with the epsilon algorithm have been
published. Thus, a detailed bibliography would be beyond the scope of
this article. Moreover, the epsilon algorithm is not restricted to
scalar sequences but can be generalized to cover for example vector
sequences. A very recent review can be found in
\cite{GMRobSal2000}.

Aitken's iterated $\Delta^2$ process (\ref{It_Aitken}) as well as Wynn's
epsilon algorithm (\ref{eps_al}) are powerful accelerators for sequences
which according to (\ref{LinLogConv}) converge linearly, and they are
also able to sum many alternating divergent series. However, they fail
completely in the case of logarithmic convergence (compare for example
\cite[Theorem 12]{Wy1966}). Moreover, in the case of divergent power
series whose series coefficients grow more strongly than factorially,
Pad\'{e} approximants or equivalently Wynn's epsilon algorithm either
converge too slowly to be numerically useful \cite{CizVrs1982,Sim1982}
or are not at all able to accomplish a summation to a unique finite
generalized limit \cite{GraGre1978}.

Brezinski showed that the inability of Wynn's epsilon algorithm of
accelerating logarithmic convergence can be overcome by a suitable
modification of the recursive scheme (\ref{eps_al}). This leads to the
so-called theta algorithm
\cite{Bre1971}:
\begin{subequations}
\label{thet_al}
\begin{eqnarray}
\vartheta_{-1}^{(n)} & = & 0 \, , \qquad
\vartheta_0^{(n)} \; = \; s_n \, , \qquad n \in \mathbb{N}_0 \, ,
\\
\vartheta_{2 k + 1}^{(n)} & = & \vartheta_{2 k-1}^{(n+1)}
\, + \, 1 / \bigl[\Delta \vartheta_{2 k}^{(n)}\bigr] \, ,
\qquad k, n \in \mathbb{N}_0 \, ,
\label{thet_al_b}
\\
\vartheta_{2 k+2}^{(n)} & = & \vartheta_{2 k}^{(n+1)} \, +
\, \frac
{\bigl[ \Delta \vartheta_{2 k}^{(n+1)} \bigr] \,
\bigl[\Delta \vartheta_{2 k + 1}^{(n+1)} \bigr]}
{\Delta^2 \vartheta_{2 k+1}^{(n)}} \, ,
\qquad k, n \in \mathbb{N}_0 \, .
\end{eqnarray}
\end{subequations}   
As in the case of Wynn's epsilon algorithm (\ref{eps_al}), only the
elements $\vartheta_{2k}^{(n)}$ with even subscripts provide
approximations to the (generalized) limit of the sequence to be
transformed.

The theta algorithm was derived with the intention of overcoming the
inability of the epsilon algorithm to accelerate logarithmic
convergence. In that respect, the theta algorithm was a great
success. Extensive numerical studies of Smith and Ford
\cite{SmiFor1982,SmiFor1979} showed that the theta algorithm is not 
only very powerful, but also much more versatile than the epsilon
algorithm. Like the epsilon algorithm, it is an efficient accelerator
for linear convergence and it is also able to sum many divergent
series. However, it is also able to accelerate the convergence of many
logarithmically convergent sequences and series. Further details as well
as additional references can be found in \cite[Section 2.9]{BreRZa1991}
or in \cite[Sections 10 and 11]{We1989}.

As for example discussed in \cite{We1991}, new sequence transformations
can be constructed by iterating explicit expressions for sequence
transformations with low transformation orders. This approach is also
possible in the case of the theta algorithm. A suitable closed-form
expression, which may be iterated, is \cite[Eq.\ (10.3-1)]{We1989}
\begin{equation}
\vartheta_2^{(n)} \; = \; s_{n+1} \, - \, \frac
{\bigl[\Delta s_n\bigr] \bigl[\Delta s_{n+1}\bigr]
\bigl[\Delta^2 s_{n+1}\bigr]}
{\bigl[\Delta s_{n+2}\bigr] \bigl[\Delta^2 s_n\bigr] -
\bigl[\Delta s_n\bigr] \bigl[\Delta^2 s_{n+1}\bigr]} \, ,
\qquad n \in \mathbb{N}_0 \, .
\label{theta_2_1}
\end{equation}
Its iteration yields the following nonlinear recursive scheme \cite[Eq.\
(10.3-6)]{We1989}:
\begin{subequations}
\label{thetit_1}
\begin{eqnarray}
\mathcal{J}_0^{(n)} & = & s_n \, , \qquad n \in \mathbb{N}_0 \, ,
\\   
\mathcal{J}_{k+1}^{(n)} & = &
\mathcal{J}_k^{(n+1)} \, - \, \frac
{\bigl[ \Delta \mathcal{J}_k^{(n)} \bigr]
\bigl[ \Delta \mathcal{J}_k^{(n+1)} \bigr]
\bigl[ \Delta^2 \mathcal{J}_k^{(n+1)} \bigr]}
{\bigl[ \Delta \mathcal{J}_k^{(n+2)} \bigr]
\bigl[ \Delta^2 \mathcal{J}_k^{(n)} \bigr]  -
\bigl[ \Delta \mathcal{J}_k^{(n)} \bigr]
\bigl[\Delta^2 \mathcal{J}_k^{(n+1)} \bigr]}
\, , \quad k,n \in \mathbb{N}_0 \, . \; \;
 \end{eqnarray}
\end{subequations}
The iterated transformation $\mathcal{J}_k^{(n)}$ has similar properties
as the theta algorithm from which it was derived: Both are very powerful
as well as very versatile. $\mathcal{J}_k^{(n)}$ is not only an
effective accelerator for linear convergence as well as able to sum many
divergent series, but it is also an effective accelerator for
logarithmic convergence
\cite{We1989,We1991,BhoBhaRoy1989,Sab1987,Sab1991,Sab1992,Sab1995}.

\typeout{==> Section 4}
\setcounter{equation}{0}
\section{Richardson Extrapolation, Wynn's Rho Algorithm, and Related
Transformations} 
\label{Sec:RichRho}

As long as $\rho$ in (\ref{LinLogConv}) is not too close to one, the
acceleration of linear convergence is comparatively simple. With the
help of Germain-Bonne's formal theory of convergence acceleration
\cite{GerB1973} and its extension \cite[Section 12]{We1989}, it can be
decided whether a sequence transformation is capable of accelerating
linear convergence or not. Moreover, several sequence transformations
are known that accelerate linear convergence effectively.

Logarithmic convergence leads to more challenging computational problems
than linear convergence. An example is the Dirichlet series
(\ref{ZetaSer}) for the Riemann zeta function. As already discussed in
Section \ref{Sec:ConstrSeqTran}, its convergence of can become so slow
that the evaluation of $\zeta (z)$ by successively adding up the terms
$(\nu+1)^{-z}$ is practically impossible.

There are also principal theoretical problems. Delahaye and
Germain-Bonne \cite{DelGerB1980,DelGerB1982} showed that no sequence
transformation can exist which is able to accelerate the convergence of
\emph{all} logarithmically convergent sequences. Consequently, an 
analogue of Germain-Bonne's beautiful formal theory of the acceleration
of linear convergence \cite{GerB1973} and its extension \cite[Section
12]{We1989} cannot exist, and the success of a convergence acceleration
process cannot be guaranteed unless additional information is available.

Nevertheless, many sequence transformations are known which work at
least for suitably restricted subsets of the class of logarithmically
convergent sequences. Examples are Richardson extrapolation
\cite{Ri1927}, Wynn's rho algorithm \cite{Wy1956b} and its iteration 
\cite[Section 6]{We1989}, Osada's modification of the rho algorithm 
\cite{Os1990}, and the modification of the $\Delta^2$ process by 
Bj{\o}rstad, Dahlquist, and Grosse \cite{BjoDahGro1981}. However, there
is considerable evidence that sequence transformations speed up
logarithmic convergence less efficiently than linear convergence (see,
for example, the discussion in \cite[Appendix A]{JeMoSoWe1999}).

Another disadvantage of logarithmic convergence is that serious
numerical instabilities are much more likely. A sequence transformation
accelerates convergence by extracting and utilizing information on the
index-dependence of the truncation errors from a finite set of input
data. This is normally accomplished by forming higher weighted
differences. If the input data are the partial sums of a \emph{strictly
alternating} series, the formation of higher weighted differences is
a remarkably stable process, but if the input data all have the
\emph{same sign}, numerical instabilities are quite likely. Thus, if the 
sequence to be transformed converges logarithmically, numerical
instabilities are to be expected, and it is usually not possible to
obtain results that are close to machine accuracy.

In some cases, these instability problems can be overcome with the help
of a condensation transformation due to Van Wijngaarden, which converts
input data having the same sign to the partial sums of an alternating
series, whose convergence can be accelerated more effectively. The
condensation transformation was first mentioned in \cite[pp.\ 126 -
127]{NaPhyLa1961} and only later published by Van Wijngaarden
\cite{vWij1965}. It was used by Daniel \cite{Dan1969} in combination
with the Euler transformation \cite{Eul1755}, and recently, it was
rederived by Pelzl and King \cite{PelKin1998}. Since the transformation
of a strictly alternating series by means of nonlinear sequence
transformations is a remarkably stable process, it was in this way
possible to evaluate special functions, that are defined by extremely
slowly convergent monotone series, not only relatively efficiently but
also close to machine accuracy \cite{JeMoSoWe1999}, or to perform
extensive quantum electrodynamical calculations
\cite{JenMogSof1999}. Unfortunately, the use of this combined
nonlinear-condensation transformation is not always possible: The
conversion of a monotone to an alternating series requires that terms of
the input series with large indices can be computed.

For the construction of sequence transformations, which are able to
accelerate logarithmic convergence, the standard interpolation and
extrapolation methods of numerical mathematics \cite{Dav1975,Joy1971}
are quite helpful. For that purpose, let us postulate the existence of a
function $\mathcal{S}$ of a continuous variable which coincides on a set
of discrete arguments $\{ x_n \}_{n=0}^{\infty}$ with the elements of
the sequence $\{ s_n \}_{n=0}^{\infty}$ to be transformed:
\begin{equation}
\mathcal{S} (x_n) \; = \; s_n \, , \qquad n \in \mathbb{N}_0 \, .
\end{equation}
This ansatz reduces the convergence acceleration problem to an
extrapolation problem. If a finite string $s_n, s_{n+1},
\ldots , s_{n+k}$ of $k+1$ sequence elements is known, one can construct
an approximation $\mathcal{S}_k (x)$ to $\mathcal{S} (x)$ which
satisfies the $k+1$ interpolation conditions
\begin{equation}
\mathcal{S}_k (x_{n+j}) \; = \; s_{n+j} \, ,
\qquad n \in \mathbb{N}_0 \, , \quad 0 \le j \le k \, .
\label{IntPolCond}
\end{equation}     
In the next step, the value of $\mathcal{S}_k (x)$ has to be determined
for $x \to x_{\infty}$. If this can be done, we can expect that
$\mathcal{S}_k (x_{\infty})$ will provide a better approximation to the
limit $s = s_{\infty}$ of the sequence $\{ s_n \}_{n=0}^{\infty}$ than
the last sequence element $s_{n+k}$ used for the construction
of $\mathcal{S}_k (x)$.

The most common interpolating functions are either polynomials or
rational functions. In the case of polynomial interpolation, it is
implicitly assumed that the $k$-th order approximant ${\cal S}_k (x)$ is
a polynomial of degree $k$ in $x$:
\begin{equation}
\mathcal{S}_k (x) \; = \, c_0 \, + \, c_1 x \, + \, \cdots
\, + \, c_k x^k \, .
\label{InterpolPol}
\end{equation}
For polynomials, the most natural extrapolation point is $x = 0$.
Accordingly, the interpolation points $x_n$ have to satisfy the
conditions    
\begin{subequations}
\label{x_n2zero}
\begin{eqnarray}
& x_0 > x_1 > \cdots > x_m > x_{m+1} > \cdots > 0 \, ,
\\
& {\displaystyle \lim_{n \to \infty} \; x_n \; = \; 0} \, .
\end{eqnarray}
\end{subequations}
The choice $x_{\infty}=0$ implies that the approximation to the limit is
to be identified with the constant term $c_0$ of the polynomial
(\ref{InterpolPol}).
   
Several different methods for the construction of interpolating
polynomials ${\cal S}_k (x)$ are known. Since only the constant term
$c_0$ of a polynomial $\mathcal{S}_k (x)$ has to be computed and since
in most applications it is desirable to compute simultaneously a whole
string of approximants $\mathcal{S}_0 (0), \mathcal{S}_1 (0),
\mathcal{S}_2 (0), \ldots$, the most economical choice is Neville's 
scheme \cite{Nev1934} for the recursive computation of interpolating
polynomials. If we set $x=0$, Neville's algorithm reduces to the
following 2-dimensional linear recursive scheme \cite[p.\ 6]{Bre1978}:
\begin{subequations}
\label{RichAl}
\begin{eqnarray}
\mathcal{N}_0^{(n)} & = \; & s_n \, ,
\quad n \in \mathbb{N}_0 \, ,
\\
\mathcal{N}_{k+1}^{(n)} & = &
\frac
{x_n \mathcal{N}_{k}^{(n+1)} \, - \,
x_{n+k+1} \mathcal{N}_{k}^{(n)} }
{x_n \, - \, x_{n+k+1}} \; ,
\quad k,n \in \mathbb{N}_0 \, . \; \; \;
\end{eqnarray}
\end{subequations}
In the literature, this variant of Neville's scheme is called Richardson
extrapolation \cite{Ri1927}.

There are functions that can be approximated more effectively by
rational functions than by polynomials. Consequently, at least for some
sequences $\{ s_n \}_{n=0}^{\infty}$ rational extrapolation should give
better results than polynomial extrapolation. Let us therefore assume
that $\mathcal{S}_k (x)$ can be expressed as the ratio of two polynomials
of degrees $l$ and $m$, respectively:
\begin{equation}
\mathcal{S}_k (x) \; = \; \frac
{a_0 + a_1 x + a_2 x^2 + \cdots + a_l x^l}
{b_0 + b_1 x + b_2 x^2 + \cdots + b_m x^m} \; ,
\qquad k, l, m \in \mathbb{N}_0 \, .
\label{InterpolRat}
\end{equation} 
This rational function contains $l + m + 2$ coefficients $a_0,
\ldots , a_l$ and $b_0, \ldots , b_m$. However, only $l + m +1$
coefficients are independent since they are determined only up to a
common nonvanishing factor. Usually, one requires either $b_0 = 1$ or
$b_m = 1$. Consequently, the $k+1$ interpolation conditions
(\ref{IntPolCond}) will determine the coefficients $a_0, \ldots, a_l$
and $b_0, \ldots, b_m$ provided that $k = l + m$ holds.

The extrapolation point $x_{\infty}=0$ is again an obvious choice. In
this case, the interpolation points $\{ x_n \}_{n=0}^{\infty}$ in
(\ref{InterpolRat}) have to satisfy (\ref{x_n2zero}), and the
approximation to the limit is to be identified with the ratio $a_0/b_0$
of the constant terms in (\ref{InterpolRat}).

If $l = m$ holds in (\ref{InterpolRat}), extrapolation to infinity is
also possible. In that case the interpolation points $\{ x_n
\}_{n=0}^{\infty}$ have to satisfy  
\begin{subequations}
\label{x_n2inf}
\begin{eqnarray}
& 0 < x_0 < x_1 < \cdots < x_m < x_{m+1} < \cdots \, ,
\\
& {\displaystyle \lim_{n \to \infty} \; x_n \; = \; \infty} \, ,
\end{eqnarray}
\end{subequations}
and the approximation to the limit is to be identified with the ratio
$a_l/b_l$ in (\ref{InterpolRat}).

As in the case of polynomial interpolation, several different algorithms
for the computation of rational interpolants are known. A discussion of
the relative merits of these algorithms as well as a survey of the
relevant literature can be found in \cite[Chapter III]{CuWu1987}. 

The most frequently used rational extrapolation technique is probably
Wynn's rho algorithm \cite{Wy1956b}:
\begin{subequations}
\label{RhoAl}
\begin{eqnarray}
\rho_{-1}^{(n)} & \; = \; & 0 \, ,
\qquad \rho_0^{(n)} \; = \; s_n \, , \qquad n \in \mathbb{N}_0 \, , \\
\rho_{k+1}^{(n)} & \; = \; & \rho_{k-1}^{(n+1)} \, + \,
\frac { x_{n+k+1} - x_n} {\rho_{k}^{(n+1)} - \rho_{k}^{(n)}} \, ,
\qquad k, n \in \mathbb{N}_0 \, .
\end{eqnarray}
\end{subequations}
Formally, the only difference between Wynn's epsilon algorithm
(\ref{eps_al}) and Wynn's rho algorithm is the sequence $\{ x_n
\}_{n=0}^{\infty}$ of interpolation points which have to satisfy
(\ref{x_n2inf}). As in the case of the epsilon algorithm, only the
elements $\rho_{2 k}^{(n)}$ with even subscripts provide approximations
to the limit.

In spite of their formal similarity, the epsilon and the rho algorithm
have complementary features. The epsilon algorithm is a powerful
accelerator for linear convergence and is also able to sum many
divergent alternating series, whereas the rho algorithm fails to
accelerate linear convergence and is not able to sum divergent
series. However, it is a very powerful accelerator for many
logarithmically convergent sequences.

The properties of Wynn's rho algorithm are discussed in books by
Brezinski \cite{Bre1977,Bre1978} and Wimp \cite{Wi1981}. Then, there is
an article by Osada \cite{Osa1996} discussing its convergence
properties, but otherwise, relatively little seems to be known about its
theoretical properties.

As in the case of Aitken's $\Delta^2$ formula (\ref{AitFor_1}), iterated
transformations can be constructed. For $k = 1$, we obtain from
(\ref{RhoAl}):
\begin{equation}
\rho_2^{(n)} \; = \; s_{n+1} \, + \, \frac
{(x_{n+2} - x_n) [ \Delta s_{n+1} ] [ \Delta s_n ] }
{ [ \Delta x_{n+1} ] [ \Delta s_{n} ] -
[ \Delta x_n ] [ \Delta s_{n+1} ] }
\, , \qquad n \in \mathbb{N}_0 \, .
\label{rho_2}
\end{equation}
This expression can be iterated yielding \cite[Eq.\ (6.3-3)]{We1991}
\begin{subequations}
\label{RhoIt}
\begin{eqnarray}
{\cal W}_{0}^{(n)} & = & s_n \, ,
\qquad n \in \mathbb{N}_0 \, ,
\\
{\cal W}_{k+1}^{(n)} & = & {\cal W}_{k}^{(n+1)} \,
+ \frac
{ (x_{n + 2 k + 2} - x_n ) \bigl[ \Delta {\cal W}_k^{(n+1)}
\bigr] \bigl[ \Delta {\cal W}_k^{(n)} \bigr] }
{(x_{n + 2 k + 2} - x_{n+1}) \bigl[ \Delta {\cal W}_k^{(n)} \bigr] -
(x_{n + 2 k + 1} - x_n) \bigl[ \Delta {\cal W}_k^{(n+1)} \bigr] } \, ,
\nonumber \\
& & k, n \in \mathbb{N}_0 \, .
\end{eqnarray}
\end{subequations}  
This is not the only possibility of iterating $\rho_2^{(n)}$. However,
the iterations derived by Bhowmick, Bhattacharya, and Roy
\cite{BhoBhaRo1989} are significantly less efficient than ${\cal
W}_k^{(n)}$, which has similar properties as Wynn's rho algorithm
\cite{We1989,We1991}.

The main practical problem with sequence transformations based
upon interpolation theory is that for a given sequence $\{ s_n
\}_{n=0}^{\infty}$ one has to find suitable interpolation
points $\{ x_n \}_{n=0}^{\infty}$ that produces good results. For
example, the Richardson extrapolation scheme (\ref{RichAl}) is normally
used in combination with the interpolation points $x_n = 1/(n+\beta)$
with $\beta > 0$. Then, $\mathcal{N}_k^{(n)}$ possesses a closed form
expression (see, for example, \cite[Lemma 2.1, p.\ 313]{MarSha1983} or
\cite[Eq.\ (7.3-20)]{We1989}),
\begin{equation}
\mathcal{N}_k^{(n)} \; = \; \Lambda_k^{(n)} (\beta, s_n) \; = \;
(-1)^k \, \sum_{j=0}^k \,
(-1)^j \, \frac {(\beta+n+j)^k} {j! \, (k-j)!} \, s_{n+j} \, ,
\qquad k, n \in \mathbb{N}_0 \, ,
\end{equation}
and the recursive scheme (\ref{RichAl}) assumes the following form
\cite[Eq.\ (7.3-21)]{We1989}:  
\begin{subequations}
\label{RichAlStand}
\begin{eqnarray}
\Lambda_0^{(n)} (\beta, s_n) & = & s_n \, ,
\qquad n \in \mathbb{N}_0 \, ,
\\
\Lambda_{k+1}^{(n)} (\beta, s_n) & = &
\Lambda_k^{(n+1)} (\beta, s_{n+1}) \, + \,
\frac {\beta+n} {k+1} \> \Delta \Lambda_k^{(n)} (\beta, s_n)
\, , \qquad k, n \in \mathbb{N}_0 \, .
\end{eqnarray}
\end{subequations}
Similarly, Wynn's rho algorithm (\ref{RhoAl}) and its iteration
(\ref{RhoIt}) are normally used in combination with the interpolation
points $x_n = n+1$, yielding the standard forms (see for example
\cite[Eq.\ (6.2-4)]{We1989})  
\begin{subequations}
\label{RhoAlStand}
\begin{eqnarray}
\rho_{-1}^{(n)} & \; = \; & 0 \, ,
\qquad \rho_0^{(n)} \; = \; s_n \, , \qquad n \in \mathbb{N}_0 \, , \\
\rho_{k+1}^{(n)} & \; = \; & \rho_{k-1}^{(n+1)} \, + \,
\frac {k+1} {\rho_{k}^{(n+1)} - \rho_{k}^{(n)}} \, ,
\qquad k, n \in \mathbb{N}_0 \, ,
\end{eqnarray}
\end{subequations}
and \cite[Section 6.3]{We1991}    
\begin{subequations}
\label{RhoItStand}
\begin{eqnarray}
{\cal W}_{0}^{(n)} & = & s_n \, ,
\qquad n \in \mathbb{N}_0 \, ,
\\
{\cal W}_{k+1}^{(n)} & = & {\cal W}_{k}^{(n+1)} \, -
\frac
{ (2 k + 2) \bigl[ \Delta {\cal W}_k^{(n+1)} \bigr]
\bigl[ \Delta {\cal W}_k^{(n)} \bigr] }
{ (2 k + 1) \Delta^2 {\cal W}_k^{(n)} } \, ,
\qquad k,n \in \mathbb{N}_0 \, .
\end{eqnarray}
\end{subequations}     
Many practically relevant logarithmically convergent sequences $\{ s_n
\}_{n=0}^{\infty}$ can be represented by series expansions of the
following kind:
\begin{equation}
s_n \; = \; s \, + \, (n+\beta)^{-\alpha} \,
\sum_{j=0}^{\infty} c_j / (n+\beta)^j \, ,
\qquad n \in \mathbb{N}_0 \, .
\label{ModSeqAlpha}
\end{equation}
Here, $\alpha$ is a positive decay parameter and $\beta$ is a positive
shift parameter. In \cite[Theorem 14-4]{We1989}, it was shown that the
standard form (\ref{RichAlStand}) of Richardson extrapolation
accelerates the convergence of sequences of the type of
(\ref{ModSeqAlpha}) if $\alpha$ is a positive integer, but fails if
$\alpha$ is nonintegral. This is also true for the standard form
(\ref{RhoAlStand}) of the rho algorithm \cite[Theorem 3.2]{Os1990}. In
the case of the iteration of Wynn's rho algorithm, no rigorous
theoretical result seems to be known but there is considerable empirical
evidence that it only works if $\alpha$ is a positive integer
\cite[Section 14.4]{We1989}).

If the decay parameter $\alpha$ of a sequence of the type of
(\ref{ModSeqAlpha}) is known, then Osada's variant of Wynn's rho
algorithm can be used \cite[Eq.\ (3.1)]{Os1990}:
\begin{subequations}
\label{OsRhoAl}
\begin{eqnarray}
{\bar \rho}_{-1}^{(n)} & \; = \; & 0 \, ,
\qquad {\bar \rho}_0^{(n)} \; = \; s_n \, ,
\qquad n \in \mathbb{N}_0 \, , \\
{\bar \rho}_{k+1}^{(n)} & \; = \; & {\bar \rho}_{k-1}^{(n+1)} \, + \,
\frac {k+\alpha} {{\bar \rho}_{k}^{(n+1)} - {\bar \rho}_{k}^{(n)}} \, ,
\qquad k, n \in \mathbb{N}_0 \, .
\end{eqnarray}
\end{subequations}
Osada also demonstrated that his variant accelerates the convergence of
sequences of the type of (\ref{ModSeqAlpha}) for arbitrary $\alpha > 0$,
and that the transformation error satisfies the following asymptotic
estimate
\cite[Theorem 4]{Os1990}:
\begin{equation}
{\bar \rho}_{2 k}^{(n)} \, - \, s \; = \;
\mathrm{O} \bigl( n^{-\alpha - 2k} \bigr) \, , \qquad n \to \infty \, .
\end{equation}

Osada's variant of the rho algorithm can be iterated. From
(\ref{OsRhoAl}) we obtain the following expression for ${\bar
\rho}_2^{(n)}$ in terms of $s_n$, $s_{n+1}$, and $s_{n+2}$:
\begin{equation}
{\bar \rho}_2^{(n)} \; = \;
s_{n+1} \, - \, \frac {(\alpha + 1)} {\alpha} \,
\frac {[\Delta s_n] [\Delta s_{n+1}]} {[\Delta^2 s_n]} \, ,
\qquad n \in \mathbb{N}_0 \, .
\end{equation}
If the iteration is done in such a way that $\alpha$ is increased by
two with every recursive step, we obtain the following recursive scheme
\cite[Eq.\ (2.29)]{We1991} which was originally derived by Bj{\o}rstad,
Dahlquist, and Grosse \cite[Eq.\ (2.4)]{BjoDahGro1981}:    
\begin{subequations}
\label{BDGal}
\begin{eqnarray}
{\overline {\cal W}}_{0}^{(n)} \; & = \; & s_n \, ,
\qquad n \in \mathbb{N}_0 \, , \\
{\overline {\cal W}}_{k+1}^{(n)} \; & = \; &
{\overline {\cal W}}_{k}^{(n+1)} \, - \,
\frac {(2 k + \alpha + 1)} {(2 k + \alpha)} \,
\frac
{\bigl[ \Delta {\overline {\cal W}}_k^{(n+1)} \bigr]
\bigl[ \Delta {\overline {\cal W}}_k^{(n)} \bigr]}
{\Delta^2 {\overline {\cal W}}_k^{(n)}} \, ,
\quad k, n \in \mathbb{N}_0 \, . \quad
\end{eqnarray}
\end{subequations}
Bj{\o}rstad, Dahlquist, and Grosse also showed that ${\overline {\cal
W}}_k^{(n)}$ accelerates the convergence of sequences of the type of
(\ref{ModSeqAlpha}), and that the transformation error satisfies the
following asymptotic estimate \cite[Eq.\ (3.1)]{BjoDahGro1981}
\begin{equation}
{\overline {\cal W}}_k^{(n)} \, - \, s \; = \;
\mathrm{O} \bigl( n^{-\alpha - 2k} \bigr) \, , \qquad n \to \infty \, .
\end{equation}  

The explicit knowledge of the decay parameter $\alpha$ is crucial for an
application of the transformations (\ref{OsRhoAl}) and (\ref{BDGal}) to
a sequence of the type of (\ref{ModSeqAlpha}).  An approximation to
$\alpha$ can be obtained with the help of the following nonlinear
transformation, which was first derived in a somewhat disguised form by
Drummond \cite{Dru1976} and later rederived by Bj{\o}rstad, Dahlquist,
and Grosse \cite{BjoDahGro1981}:
\begin{equation}
T_n \; = \; \frac
{ [\Delta^2 s_n ] \, [\Delta^2 s_{n+1} ]}
{[\Delta s_{n+1} ] \, [\Delta^2 s_{n+1} ] \, - \,
[\Delta s_{n+2} ] \, [\Delta^2 s_{n}]} \; - \; 1 \, ,
\qquad n \in \mathbb{N}_0 \, .
\label{DecPar}
\end{equation}
$T_n$ is essentially a weighted $\Delta^3$ method, which implies that it
is potentially very unstable. Thus, stability problems are likely to
occur if the relative accuracy of the input data is low. Bj{\o}rstad,
Dahlquist, and Grosse \cite[Eq.\ (4.1)]{BjoDahGro1981} also showed that
\begin{equation}
\alpha \; = \; T_n \, + \, O (1/n^2) \,,\qquad n \to \infty \, ,
\label{T_n}
\end{equation}
if the elements of a sequence of the type of (\ref{ModSeqAlpha}) are
used as input data.

\typeout{==> Section 5}
\setcounter{equation}{0}
\section{Transformations with Explicit Remainder Estimates} 
\label{Sec:LevinType}

As discussed before, the action of a sequence transformation corresponds
at least conceptually to the construction an approximations to the
truncation error $r_n$, whose elimination from $s_n$ leads to an
acceleration of convergence or a summation. However, the remainders
$r_n$ may depend on $n$ in a very complicated way, and the construction
of approximations to the $r_n$ and their subsequent elimination can be
very difficult. In some cases, however, structural information on the
$n$-dependence of the $r_n$ is available. For example, the truncation
error of a convergent series with strictly alternating and monotonously
decreasing terms is bounded in magnitude by the first term not included
in the partial sum and has the same sign as this term
\cite[p.\ 132]{Kno1964}. The first term neglected is also the best
simple estimate for the truncation error of a strictly alternating
nonterminating hypergeometric series ${}_2 F_0 (\alpha, \beta; - x)$
with $\alpha, \beta, x > 0$ \cite[Theorem 5.12-5]{Car1977}, which
diverges for every nonzero argument $x$. Such an information should be
extremely helpful. Unfortunately, the sequence transformations
considered so far cannot benefit from it.

Structural information of that kind can be incorporated into the
transformation process via explicit remainder estimates $\{ \omega_n
\}_{n=0}^{\infty}$, as it was first done by Levin \cite{Lev1973}. For
that purpose, let us assume that the remainders $r_n$ of a sequence $\{
s_n \}_{n=0}^{\infty}$ can be partitioned into a \emph{remainder
estimate} $\omega_n$ multiplied by a \emph{correction term} $z_n$
according to $r_n = \omega_n z_n$. The remainder estimates are chosen
according to some rule and may depend on $n$ in a very complicated
way. If the remainder estimates $\{ \omega_n\}_{n=0}^{\infty}$ correctly
describe the essential features of the remainders $\{ r_n
\}_{n=0}^{\infty}$, then the $\{ z_n \}_{n=0}^{\infty}$ should depend on 
$n$ in a relatively smooth way. Of course, we tacitly assume here that
the products $\omega_n z_n$ are in principle capable of producing
sufficiently accurate approximations to the remainders.

Thus, we have to find a sequence transformation that is exact for the
model sequence
\begin{equation}
\tilde{s}_n \; = \; \tilde{s} \, + \, \omega_n \, z_n \, ,
\qquad n \in \mathbb{N}_0 \, ,
\label{Mod_Seq_Om}
\end{equation}
where the remainder estimates $\{ \omega_n \}_{n=0}^{\infty}$ are
assumed to be known. The principal advantage of this approach is that
only approximations to the correction terms $\{ z_n \}_{n=0}^{\infty}$
have to be determined. If good remainder estimates can be found, the
determination of $z_n$ and the subsequent elimination of $\omega_n z_n$
from $s_n$ often leads to clearly better results than the construction
and elimination of other approximations to $r_n$.  The explicit
utilization of information contained in remainder estimates is the major
difference between the sequence transformations discussed in this
Section and the other transformations of this article.
 
The model sequence (\ref{Mod_Seq_Om}) has another indisputable
advantage: There is a \emph{systematic} way of constructing a sequence
transformation which is exact for this model sequence. Let us assume
that a \emph{linear} operator ${\hat T}$ can be found which annihilates
the correction term $z_n$ according to ${\hat T} (z_n) = 0$. Then, a
sequence transformation, which is exact for the model sequence
(\ref{Mod_Seq_Om}), can be obtained by applying ${\hat T}$ to
$[\tilde{s}_n - \tilde{s}] / \omega_n = z_n$. Since ${\hat T}$
annihilates $z_n$ and is by assumption linear, the following sequence
transformation ${\cal T}$ is \emph{exact} for the model sequence
(\ref{Mod_Seq_Om}) \cite[Eq.\ (3.2-11)]{We1989}:
\begin{equation}
{\cal T} (\tilde{s}_n, \omega_n) \; = \; \frac
{{\hat T} (\tilde{s}_n / \omega_n )} {{\hat T} (1 / \omega_n )}
\; = \; \tilde{s} \, .
\label{GenSeqTr}
\end{equation}
Originally, the construction of sequence transformations via
annihilation operators was introduced in \cite[Section 3.2]{We1989} in
connection with a rederivation of Levin's transformation
\cite{Lev1973} and the construction of some other, closely related sequence
transformations \cite[Sections 7 - 9]{We1989}. Later, this operator
approach was also used and discussed by Brezinski
\cite{Bre2000,Bre1996,Bre1997}, Brezinski and Redivo Zaglia 
\cite{BreRZa1991,BreRZa1994a,BreRZa1994b}, Brezinski and
Salam \cite{BreSal1995}, Brezinski and Matos
\cite{BreMa1996}, Matos \cite{Mat2000}, and Homeier \cite{Hom2000}.

If the annihilation operator ${\hat T}$ in (\ref{GenSeqTr}) is based
upon the finite difference operator $\Delta$, simple and yet very
powerful sequence transformations are obtained
\cite[Sections 7 - 9]{We1989}. As is well known, $\Delta^k$ annihilates
a polynomial $P_{k-1} (n)$ of degree $k - 1$ in $n$. Thus, the
correction terms should be chosen in such a way that multiplication of
$z_n$ by some $w_k (n)$ yields a polynomial $P_{k-1} (n)$ of degree
$k-1$ in $n$.

Since $\Delta^k w_k (n) z_n = \Delta^k P_{k-1} (n) = 0$, the weighted
difference operator ${\hat T} = \Delta^k w_k (n)$ annihilates $z_n$, and
the corresponding sequence transformation (\ref{GenSeqTr}) is given by
the ratio
\begin{equation}
{\cal T}_k^{(n)} \bigl( w_k (n) \big\vert s_n, \omega_n \bigr)
\; = \; \frac
{\Delta^k \{ w_k (n) s_n / \omega_n \}}
{\Delta^k \{ w_k (n) / \omega_n \}} \, .
\label{Seq_w_k}
\end{equation}   
Several different sequence transformations are obtained by specializing
$w_k (n)$. For instance, $w_k (n) = (n + \zeta)^{k-1}$ with $\zeta > 0$
yields Levin's sequence transformation \cite{Lev1973}:
\begin{eqnarray}
\lefteqn{
{\cal L}_{k}^{(n)} \bigl(\zeta, s_n, \omega_n\bigr) 
\; = \; \frac
{ \Delta^k \, \{ (n + \zeta)^{k-1} \> s_n / \omega_n\} }
{ \Delta^k \, \{ (n + \zeta)^{k-1}  / \omega_n \} }
} \nonumber \\
& \; = \; \frac
{\displaystyle
\sum_{j=0}^{k} \; ( - 1)^{j} \; {{k} \choose {j}} \;
\frac
{(\zeta + n +j )^{k-1}} {(\zeta + n + k )^{k-1}} \;
\frac {s_{n+j}} {\omega_{n+j}} }
{\displaystyle
\sum_{j=0}^{k} \; ( - 1)^{j} \; {{k} \choose {j}} \;
\frac
{(\zeta + n +j )^{k-1}} {(\zeta + n + k )^{k-1}} \;
\frac {1} {\omega_{n+j}} } \, .
\label{LevTr}
\end{eqnarray}   
The shift parameter $\zeta$ has to be positive to allow $n = 0$ in
(\ref{LevTr}). The most obvious choice is $\zeta = 1$. According to
Smith and Ford \cite{SmiFor1982,SmiFor1979}, Levin's transformation is
among the most powerful and most versatile sequence transformations that
are currently known.

The unspecified weights $w_k (n)$ in (\ref{Seq_w_k}) can also be
chosen to be Pochhammer symbols according to $w_k (n) = (n +
\zeta)_{k-1}$ with $\zeta > 0$, yielding \cite[Eq.\ (8.2-7)]{We1989} 
\begin{eqnarray}
\lefteqn{
{\cal S}_{k}^{(n)} \bigl(\zeta , s_n, \omega_n\bigr) \; = \; \frac
{ \Delta^k \, \{ (n + \zeta)_{k-1} \> s_n / \omega_n\} }
{ \Delta^k \, \{ (n + \zeta)_{k-1}  / \omega_n \} }
} \nonumber \\
& \; = \; \frac
{\displaystyle
\sum_{j=0}^{k} \; ( - 1)^{j} \; {{k} \choose {j}} \;
\frac {(\zeta + n +j )_{k-1}} {(\zeta + n + k )_{k-1}} \;
\frac {s_{n+j}} {\omega_{n+j}} }
{\displaystyle
\sum_{j=0}^{k} \; ( - 1)^{j} \; {{k} \choose {j}} \;
\frac {(\zeta + n +j )_{k-1}} {(\zeta + n + k )_{k-1}} \;
\frac {1} {\omega_{n+j}} } \, .
\label{SidTr}
\end{eqnarray}
As shown in several articles, this transformation is very effective if
strongly divergent alternating series are to be summed
\cite{We1989,We1994a,We1990,WeCi1990,We1996c,JeMoSoWe1999,WeCiVi1991,%
CiViWe1991,We1992,WeCiVi1993,CiViWe1993,We1996a,We1996b,We1996d,%
We1997,JeBeWeSo2000,RoyBhaBho1996,BhaBhoRoy1997,RoyBhaBho1998}. Again,
the most obvious choice for the shift parameter is $\zeta = 1$.

Levin's transformation ${\cal L}_{k}^{(n)} \bigl(\zeta, s_n,
\omega_n\bigr)$ is by construction exact for the model sequence $s_n = s 
+ \omega_n \sum_{j=0}^{k-1} c_j/(n + \zeta)^j$. Consequently, Levin's
transformation should work well if the ratio $[s_n - s]/\omega_n$ can be
expressed as a power series in $1/(n + \zeta)$.  Similarly, ${\cal
S}_{k}^{(n)} \bigl(\zeta , s_n, \omega_n \bigr)$ is by construction
exact for the model sequence $s_n = s + \omega_n \sum_{j=0}^{k-1} c_j /
(n + \zeta)_j$, which is a truncated factorial series \cite[Eq.\
(8.2-1)]{We1989}. Accordingly, ${\cal S}_{k}^{(n)} \bigl(\zeta, s_n,
\omega_n\bigr)$ should give good results if the ratio $[s_n -
s]/\omega_n$ can be expressed as a factorial series.

Given a suitable sequence $\{ \omega_n \}_{n=0}^{\infty}$ of remainder
estimates, the sequence transformations ${\cal L}_{k}^{(n)}$ and ${\cal
S}_{k}^{(n)}$ can be computed via their explicit expressions
(\ref{LevTr}) and (\ref{SidTr}). However, it is usually more effective
to compute the numerator and denominator sums in (\ref{LevTr}) and
(\ref{SidTr}) with the help of three-term recursions \cite[Eqs.\ (7.2-8)
and (8.3-7)]{We1989}.

The explicit incorporation of the information contained in the remainder
estimates makes the transformations (\ref{LevTr}) and (\ref{SidTr})
potentially very powerful. However, this is also their major potential
weakness. If remainder estimates can be found such that the products
$\omega_n z_n$ provide good approximations to the remainders, Levin-type
sequence transformations should work very well. If, however, good
remainder estimates cannot be found, sequence transformations of that
kind perform poorly. Consequently, a fortunate choice of the remainder
estimates in (\ref{LevTr}) and (\ref{SidTr}) is of utmost importance
since it ultimately determines success or failure.

The difference operator $\Delta$ is linear. Consequently, the effect of
the general sequence transformation (\ref{Seq_w_k}) on an arbitrary
sequence $\{ s_n \}_{n=0}^{\infty}$ with (generalized) limit $s$ can be
expressed as follows:
\begin{equation}
{\cal T}_k^{(n)} \bigl( w_k (n) \big\vert s_n, \omega_n \bigr) \; = \;
s \, + \, {\displaystyle \frac
{\Delta^k \{ w_k (n) (s_n - s) / \omega_n \}}
{\Delta^k \{ w_k (n) / \omega_n \}}} \, .
\end{equation}
Obviously, ${\cal T}_k^{(n)} \bigl( w_k (n) \big\vert s_n, \omega_n
\bigr)$ converges to $s$ if the ratio on the right-hand side can be made
arbitrarily small. This is the case if $\Delta^k w_k (n)$ annihilates
$[s_n - s]/\omega_n$ more effectively than $1/\omega_n$. Thus, one
should try to find remainder estimates such that the ratios $[s_n -
s]/\omega_n$ depend on $n$ only weakly:
\begin{equation}
[s_n - s] / \omega_n \; = \;  c + O (1/n) \, , \quad c \ne 0 \, ,
\quad n \to \infty \, .
\label{asy_Rem_Est}
\end{equation}    
This asymptotic condition does not determine the remainder estimates
uniquely. Therefore, it is at least in principle possible to find for a
given sequence an unlimited variety of different remainder estimates,
which all satisfy this asymptotic condition.

On the basis of heuristic and asymptotic arguments, Levin \cite{Lev1973}
suggested the following simple remainder estimates, which nevertheless
often work remarkably well:
\begin{eqnarray}
\omega_n & \; = \; & (\zeta + n) \, \Delta s_{n-1} \, , 
\qquad \zeta > 0 \, ,
\label{u_Est} \\
\omega_n & \; = \; & \Delta s_{n-1} \, , \label{t_Est} \\
\omega_n & \; = \; &  \frac
{\Delta s_{n-1} \Delta s_n} {\Delta s_{n-1} - \Delta s_n} \, .
\label{v_Est}
\end{eqnarray}  
The use of these remainder estimates in (\ref{LevTr}) yields Levin's
$u$, $t$, and $v$ transformation, respectively
\cite[Eqs.\ (7.3-5), (7.3-7), and (7.3-11)]{We1989}:
\begin{eqnarray}
u_k^{(n)} \bigl(\zeta, s_n\bigr) & = &
{\cal L}_{k}^{(n)} 
\bigl(\zeta, s_n, (\zeta + n) \Delta s_{n-1}\bigr) \, ,
\\
t_k^{(n)} \bigl(\zeta, s_n\bigr) & = &
{\cal L}_{k}^{(n)} \bigl(\zeta, s_n, \Delta s_{n-1}\bigr) \, ,
\\
v_k^{(n)} \bigl(\zeta, s_n\bigr) & = &
{\cal L}_{k}^{(n)} (\zeta, s_n, \bigl(\Delta s_{n-1} \Delta s_n\bigr) /
(\Delta s_{n-1} - \Delta s_n)) \, .
\end{eqnarray}   
Later, Smith and Ford \cite{SmiFor1979} suggested the remainder estimate
\begin{equation}
\omega_n \; = \; \Delta s_n \, ,
\label{d_Est}
\end{equation}
which yields Levin's $d$ transformation \cite[Eq.\ (7.3-9)]{We1989}:
\begin{equation}
d_k^{(n)} \bigl(\zeta, s_n\bigr) \; = \;
{\cal L}_{k}^{(n)} \bigl(\zeta, s_n, \Delta s_n\bigr) \, .
\label{dLevTr}
\end{equation}   
Levin's $t$ and $d$ transformations are capable of accelerating linear
convergence and they are particularly efficient in the case of
alternating series, but fail to accelerate logarithmic
convergence. Levin's $u$ and $v$ transformation are more versatile since
they not only accelerate linear convergence but also many
logarithmically convergence sequences and series. A more detailed
discussion of the properties of these remainder estimates, some
generalizations, additional heuristic motivation, and a description of
the types of sequences, for which these estimates should be effective,
can be found in \cite[Sections 7 and 12 - 14]{We1989}. The main
advantage of the simple remainder estimates (\ref{u_Est}) -
(\ref{v_Est}) and (\ref{d_Est}) is that they can be used in situations
in which only the numerical values of a few elements of a slowly
convergent or divergent sequence are known.

The remainder estimates (\ref{u_Est}) - (\ref{v_Est}) and (\ref{d_Est})
can also be used in the case of the sequence transformation
(\ref{SidTr}), yielding \cite[Eqs.\ (8.4-2), (8.4-3), (8.4-4), and
(8.4-5)]{We1989}
\begin{eqnarray}
y_k^{(n)} \bigl(\zeta, s_n\bigr) & = &
{\cal S}_{k}^{(n)} 
\bigl(\zeta, s_n, (\zeta + n) \Delta s_{n-1}\bigr) \, ,
\\
{\tau}_k^{(n)} \bigl(\zeta, s_n\bigr) & = &
{\cal S}_{k}^{(n)} \bigl(\zeta, s_n, \Delta s_{n-1}\bigr) \, ,
\\
{\phi}_k^{(n)} \bigl(\zeta, s_n\bigr) & = &
{\cal S}_{k}^{(n)} \bigl(\zeta, s_n, (\Delta s_{n-1} \Delta s_n) /
(\Delta s_{n-1} - \Delta s_n)\bigr) \, ,
\\
{\delta}_k^{(n)} \bigl(\zeta, s_n\bigr) & = &
{\cal S}_{k}^{(n)} \bigl(\zeta, s_n, \Delta s_n\bigr) \, .
\label{dSidTr}
\end{eqnarray} 
Alternative remainder estimates for the sequence transformations
(\ref{LevTr}) and (\ref{SidTr}) are discussed in
\cite{We1989,HoWe1995}. Convergence properties of the sequence 
transformations (\ref{LevTr}) and (\ref{SidTr}) and their variants were
analyzed in articles by Sidi \cite{Sid1979,Sid1980a,Sid1986}, in
\cite[Sections 12 - 14]{We1989}, in \cite[Section 4]{WeCiVi1993}, and
also in \cite{Hom2000}.               

Other sequence transformations, which are also special cases of the
general sequence transformation (\ref{GenSeqTr}), can be found in
\cite[Sections 7 - 9]{We1989}, in the book by Brezinski and Redivo
Zaglia \cite[Section 2.7]{BreRZa1991}, or in a recent review by Homeier
\cite{Hom2000}, which is the currently most complete source of
information on Levin-type transformations. A sequence transformation,
which interpolates between the transformations (\ref{LevTr}) and
(\ref{SidTr}), was described in \cite{We1992}.

\typeout{==> Section 6}
\setcounter{equation}{0}
\section{Numerical Aspects}
\label{Sec:NumAsp}

As discussed before, sequence transformations try to accomplish an
acceleration of convergence or a summation by detecting and utilizing
regularities in the behavior of the elements of the sequence to be
transformed. For sufficiently large indices $n$, one can expect that
certain asymptotic regularities do exist. However, sequence
transformations are normally used with the intention of avoiding the
asymptotic domain, i.e., one tries to construct the transforms from the
\emph{leading} elements of the input sequence. Unfortunately, sequence
elements $s_n$ with small indices $n$ often behave
irregularly. Consequently, it can happen that a straightforward
application of a sequence transformation is ineffective and even leads
to completely nonsensical results. In fact, one should not be too
surprised that a strategy, which tries to avoid the asymptotic domain by
extracting asymptotic information from the leading elements of an input
sequence, occasionally runs into trouble.

As a possible remedy, one should analyze the behavior of the input data
as a function of the index and exclude highly irregular sequence
elements -- usually the leading ones -- from the transformation
process. This gives a much better chance of obtaining good and reliable
transformation results. 

In view of the fact that irregular input data can never be excluded, one
might expect that there are many references dealing with this
topic. However, apart from a recent reference of my own \cite{We2001}, I
am only aware of an article by Gander, Golub, and Gruntz
\cite{GaGoGr1990} where it is shown that the convergence of
extrapolations of iteration sequences can be improved by excluding the
leading elements of the input sequence from the extrapolation
process. Otherwise, a discussion of the impact of irregular input data
on convergence acceleration and summation processes seems to be part of
the oral tradition only.

A numerical process can only involve a finite number of arithmetic
operations. Thus, a sequence transformation $\mathcal{T}$ can only use
finite subsets of the type $\{ s_n, s_{n+1}, \ldots, s_{n+l}
\}$ for the computation of a new sequence element
$s^{\prime}_n$. Moreover, all the commonly used sequence transformations
$\mathcal{T}$ can be represented by infinite sets of doubly indexed
quantities $T_{k}^{(n)}$ with $k, n \in \mathbb{N}_0$ that can be
displayed in a two-dimensional array called the \emph{table} of
$\mathcal{T}$. The superscript $n$ always indicates the minimal index
occurring in the finite subset of input data used for the computation of
$T_k^{(n)}$. The subscript $k$ -- usually called the \emph{order} of the
transformation -- is a measure for the complexity of the transformation
process which yields $T_k^{(n)}$. The elements $T_k^{(n)}$ are gauged in
such a way that $T_0^{(n)}$ corresponds to an untransformed sequence
element according to $T_0^{(n)} = s_n$. An increasing value of $k$
implies that the complexity of the transformation process
increases. Moreover, $l = l (k)$ in the substring $\{ s_n, s_{n+1},
\ldots, s_{n+l} \}$ also increases. This means that for every 
$k, n \in \mathbb{N}_0$ the sequence transformation $\mathcal{T}$
produces a new transform according to
\begin{equation}
 T_k^{(n)} \; = \;
{\cal T} \bigl(s_n, s_{n+1}, \ldots , s_{n+l(k)} \bigr) \, .
\end{equation}
The relationship, which connects $k$ and $l$, is specific for a given
sequence transformation ${\cal T}$.

Let us assume that a sequence transformation $\mathcal{T}$ is to be used
to speed up the convergence of a sequence $\{ s_n \}_{n=0}^{\infty}$ to
its limit $s = s_{\infty}$. One can try to obtain a better approximation
to $s$ by proceeding on an unlimited variety of different \emph{paths}
in the table of $\mathcal{T}$. Two extreme types of paths -- and also
the ones which are predominantly used in practical applications -- are
\emph{order-constant} paths $\bigl\{ T_{k}^{(n+\nu)}
\bigr\}_{\nu=0}^{\infty}$ with fixed transformation order $k$ and
increasing superscript, and \emph{index-constant} paths $\bigl\{
T_{k+\kappa}^{(n)} \bigr\}_{\kappa=0}^{\infty}$ with fixed minimal index
$n$ and increasing subscript.

These two types of paths differ significantly. In the case of an
order-constant path, a \emph{fixed} number of $l+1$ sequence elements
$\{ s_n, s_{n+1}, \ldots s_{n+l} \}$ is used for the computation of
$T_k^{(n)}$, and the starting index $n$ of this string of fixed length
is increased successively until either convergence is achieved or the
available elements of the input sequence are used up. In the case of
an index-constant path, the starting index $n$ is kept fixed at a low
value (usually $n = 0$ or $n = 1$) and the transformation order $k$ is
increased and with it the number of elements contained in the subset $\{
s_n, s_{n+1}, \ldots s_{n+l(k)} \}$. Thus, on an index-constant path
$T_k^{(n)}$ is always computed with the highest possible transformation
order $k$ from a given set of input data.

In order to clarify the differences between order-constant and
index-constant paths, let us consider the computation of Pad\'e
approximants with the help of Wynn's epsilon algorithm
(\ref{eps_al}). It follows from (\ref{Eps_Pade}) that the epsilon
algorithm (\ref{eps_al}) effects the following transformation of the
partial sums (\ref{Par_Sum_PS}) of the power series for some function
$f$ to Pad\'{e} approximants:
\begin{equation}
\bigl\{ f_{n} (z), f_{n+1} (z), \ldots, f_{n+2k} (z) \bigr\}
\; \longrightarrow \; [ n + k / k ] \, .
\end{equation}
Thus, if we use a window consisting of $2k + 1$ partial sums $f_{n+j}
(z)$ with $0 \le j \le 2k$ on an order-constant path and increase the
minimal index $n$ successively, we obtain the following sequence of
Pad\'{e} approximants:
\begin{equation}
[ n + k / k ], [ n + k + 1 / k ], \ldots, [ n + k + m / k ],
\ldots \, .
\end{equation}
Only $2k + 1$ partial sums are used for the computation of the Pad\'{e}
approximants, although many more may be known. Obviously, the available
information is not completely utilized on such an order-constant path.

Then, the degree of the numerator polynomial of $[ n + k + m / k ]$
increases with increasing $m \in \mathbb{N}_0$, whereas the degree of
the denominator polynomial remains fixed. Thus, these Pad\'{e}
approximants look unbalanced. Instead, it seems to be much more natural
to use \emph{diagonal} Pad\'{e} approximants, i.e., Pad\'{e}
approximants with numerator and denominator polynomials of equal degree,
or -- if this is not possible -- to use Pad\'{e} approximants with
degrees of the numerator and denominator polynomials that differ as
little as possible.

This approach has many theoretical as well as practical advantages. Wynn
could show that if the partial sums $f_0 (z)$, $f_1 (z)$, $\cdots$,
$f_{2 n} (z)$ of a Stieltjes series are used for the computation of
Pad\'{e} approximants, then the diagonal approximant $[ n / n ]$
provides the most accurate approximation to the corresponding Stieltjes
function $f (z)$, and if the partial sums $f_0 (z)$, $f_1 (z)$,
$\cdots$, $f_{2 n + 1} (z)$ are used for the computation of Pad\'{e}
approximants, then for $z > 0$ either $[ n + 1 / n ]$ or $[ n / n + 1 ]$
provides the most accurate approximation \cite{Wy1968}. A detailed
discussion of Stieltjes series and their special role in the theory of
Pad\'{e} approximants can be found in \cite[Section 5]{BaGM1996}.

Thus, it is apparently an obvious idea to try to use either diagonal
Pad\'{e} approximants or their closest neighbors whenever possible. If
the partial sums $f_0 (z)$, $f_1 (z)$, $\ldots$, $f_m (z)$, $\ldots$ are
computed successively and used as input data in the epsilon algorithm
(\ref{eps_al}), then we obtain the following staircase sequence in the
Pad\'{e} table \cite[Eq.\ (4.3-7)]{We1989}, which exploits the available 
information optimally:
\begin{equation}
[0/0], [1/0], [1/1], \ldots 
[\nu / \nu], [\nu + 1/ \nu], [\nu +1/ \nu +1], \ldots \, .
\end{equation}
This example indicates that index-constant paths are in principle
computationally more efficient than order-constant paths since they
exploit the available information optimally. This is also true for all
other sequence transformations considered in this article.

Another serious disadvantage of order-constant paths is that they cannot
be used for the summation of divergent sequences and series since
increasing $n$ in the set $\{ s_n, s_{n+1}, \ldots, s_{n+l} \}$ of input
data normally only increases divergence. Thus, it is apparently an
obvious idea to use exclusively index-constant paths, and preferably
those which start at a very low index $n$, for instance at $n = 0$ or $n
= 1$. This is certainly a good idea if all elements of the input
sequence contain roughly the same amount of useful information. If,
however, the leading terms of the sequence to be transformed behave
irregularly, they cannot contribute useful information, or -- to make
things worse -- they contribute \emph{wrong} information. In such a case
it is usually necessary to exclude the leading elements of the input
sequence from the transformation process. Then, one should use either an
order-constant path or an index-constant path with a sufficiently large
starting index $n$. The use of an order-constant path has the additional
advantage that the diminishing influence of irregular input data with
small indices $n$ should become obvious from the transformation results
as $n$ increases.

In \cite{WeCiVi1991,WeCiVi1993,We1996b,We1996d} I did extensive quantum
mechanical calculations for anharmonic oscillators by summing strongly
divergent perturbation expansions. The coefficients of these
perturbation expansions are all rational numbers. Thus, they can be
computed free of errors with the help of the exact rational arithmetics
of a computer algebra system like Maple. In such a case, it is not only
an obvious idea but in fact necessary to do the summation calculations
on an index-constant path. The maximum transformation order in these
summation calculations was only limited by the number of perturbative
coefficients that could be computed before the available memory of the
computer was exhausted.

The situation was completely different when I was involved in the
extrapolation of quantum chemical oligomer calculations to the infinite
chain limit \cite{WeLie1990,CioWe1993,WeKi2000}. As discussed before,
sequence transformations have to access the information stored in the
later digits of the input data. However, the input data are produced by
molecular \emph{ab initio} programs which are huge packages of FORTRAN
code (more than $100~000$ lines of code). Normally, these programs
operate in DOUBLE PRECISION which corresponds to an accuracy of 14 - 16
decimal digits (depending on the compiler). Unfortunately, this does not
mean that their results have this accuracy. Some of the leading digits
may have converged, and some of the trailing digits are corrupted
because of internal approximations used in the molecular program or
because of inevitable rounding errors. Thus, in most cases there is only
a relatively narrow window of digits that can provide useful information
for the transformation process.  Moreover, the complexity of the
calculations done in such a molecular program makes it impossible to
obtain realistic estimates of the errors by the standard approaches of
numerical mathematics. In such a case, it is recommendable to do the
extrapolations on order-constant paths with low transformation orders,
because high transformation orders can easily lead to meaningless
extrapolation results. Further details can be found in \cite{WeKi2000}.

It may look like an obvious idea to use the results of a transformation
process as input data for another sequence transformation. At least from
a conceptual point of view, this is essentially the same as iterating a
low order transformation like Aitken's $\Delta^2$ formula
(\ref{AitFor_1}). If, however, the input data for the second
transformation were obtained on an index-constant path by the first
transformation, then I see basically two problems. The law, which
governs the convergence of the first transformation, may be unknown
and/or hard to find. Thus, it may be difficult to find a suitable second
transformation. The second and -- as I feel -- more serious restriction
is numerical in nature. The trailing digits of the input data for the
second transformation may be largely corrupted by the first
transformation. Thus, they cannot provide the information needed for a
further improvement of convergence. There are articles which describe
successful attempts of combining different sequence transformations
\cite{JeWeS02000,Bel1989}. However, I suspect that many more attempts
were not reported in the literature because they failed to accomplish
something substantial.

\typeout{==> Section 7}
\setcounter{equation}{0}
\section{Outlook}
\label{Sec:Outlook}

If a review article is to written and if there is only a limited amount
of space available, there is always the danger that the author
emphasizes those aspects he knows particularly well, whereas other
aspects are not treated as thoroughly as they probably should. This is
of course also true for this review. So, I will now try to give a more
balanced view by mentioning some additional applications of sequence
transformations.

In this review, exclusively the transformation of sequences of real or
complex numbers was treated. However, sequence transformations can also
be formulated for vector or matrix problems. Sequence transformations
for vector sequences are for instance treated in the books by Brezinski
and Redivo Zaglia \cite[Section 4]{BreRZa1991} and Brezinski
\cite{Bre1997}, or in articles by MacLeod \cite{MacL1986}, 
Graves-Morris and Saff \cite{GraMoSaf1988,GraMoSaf1991}, Sidi
\cite{Sid1988/9}, Smith, Ford, and Sidi 
\cite{SmiForSid1986,SmiForSid1988}, Sidi and Bridger \cite{SidBri1988},  
Jbilou and Sadok \cite{JbiSad1991,JbiSad1995},
Osada \cite{Osa1991,Osa1992,Osa1996a}, Brezinski and Sadok 
\cite{BreSad1992}, Matos \cite{Mat1992}, Graves-Morris 
\cite{GraMo1992,GraMo1996}, Graves-Morris and Roberts 
\cite{GraMoRob1994,GraMoRob1997}, Brezinski and Salam \cite{BreSal1995},
Homeier, Rast, and Krienke \cite{HomRasKri1995}, Salam
\cite{Sal1996,Sal1998,Sal1999}, Graves-Morris and Van Iseghem
\cite{GraMoVIs1997}, Roberts \cite{Rob1998a,Rob1998b}, and Graves-Morris,
Roberts, and Salam \cite{GMRobSal2000},

In the practice of the author, convergence acceleration and summation of
infinite series has dominated. However, sequence transformations can
also be very useful in the context of numerical quadrature, in particular
if an oscillatory function is to be integrated over a semiinfinte
interval. The use of sequence transformations in connection with
numerical quadratures is for instance discussed in a book by Evans
\cite{Eva1993}, or in articles by Sidi
\cite{Sid1980b,Sid1980c,Sid1982a,Sid1982b,Sid1987,Sid1988,Sid1990,Sid1999},  
Levin and Sidi \cite{LevSid1981}, Levin \cite{Lev1982}, Greif and Levin
\cite{GreLev1998}, Safouhi and Hoggan
\cite{SafHog1998,SafHog1999a,SafHog1999b,SafHog1999c}, Safouhi, Pinchon,
and Hoggan \cite{SafPinHog1998}, and Safouhi \cite{Saf2000,Saf2001}.

So far, sequence transformations were typically used to obtain better
approximations to the limits of sequences, series, or numerical
quadrature schemes. However, it is possible to use Pad\'{e} approximants
or other sequence transformations to make predictions for unknown
perturbation series coefficients. For example, in some subfields of
theoretical physics it is extremely difficult to compute more than just
a few perturbative coefficients, and even these few coefficients are
usually affected by large relative errors. Thus, sequence
transformations -- mainly Pad\'{e} approximants -- were used to make
predictions for unknown coefficients or to check the accuracy of already
computed coefficients. Further details as well as many references can be
found in \cite{JeBeWeSo2000,JeWeS02000,We2000,ChiEliMirSte2000}. In the
majority of these predictions calculations, rational expressions in an
unspecified symbolic variable were constructed, which were then expanded
in a Taylor polynomial of suitable length with the help of a computer
algebra program like Maple or Mathematica. However, in the case of
Aitken's iterated $\Delta^2$ process (\ref{It_Aitken}), of Wynn's
epsilon algorithm (\ref{eps_al}), and of the iteration (\ref{thetit_1})
of Brezinski's theta algorithm explicit recursive schemes for
predictions were recently derived \cite{We2000}. With the help of these
recursive schemes, it is possible to make predictions for power series
coefficients with very large indices. In this way, it was possible to
gain some insight into the mathematical properties of a perturbation
expansion for a non-Hermitian $\mathcal{PT}$-symmetric Hamiltonian
\cite{BenWen2000}.

I hope that this review could show that sequence transformations are
extremely useful numerical tools, and that there is a lot of active
research going on, both on the mathematical properties of sequence
transformations as well as on their application in applied mathematics
and in the mathematical treatment of the sciences. Moreover, I am
optimistic and expect further progress in the near future.

\typeout{==> Acknowledgments}
\setcounter{equation}{0}
\section*{Acknowledgments}
\label{Sec:Ack}

I would like to thank Professor Dzevad Belki\'{c} for his invitation to
contribute to this topical issue of Journal of Computational Methods in
Sciences and Engineering. Financial support by the Fonds der Chemischen
Industrie is gratefully acknowledged.

\typeout{==> References}
\small
\parsep0.15ex plus0.04ex minus0.04ex
\itemsep0ex plus0.04ex  

\end{document}